\documentclass[a4paper,12pt]{article}
\usepackage[english]{babel}
\usepackage{tgpagella,eulervm}
\usepackage{amsmath}
\usepackage{color}
\usepackage{xcolor}

\definecolor{darkblue}{rgb}{0.0, 0.0, 0.8}
\definecolor{darkred}{rgb}{0.8, 0.0, 0.0}
\definecolor{darkgreen}{rgb}{0.6, 0.15, 0.15}

\usepackage{url}
\usepackage[normalem]{ulem}
\usepackage[utf8x]{inputenc}
\usepackage[T1]{fontenc}
\usepackage{caption}
\usepackage{hyperref}
\hypersetup{
    colorlinks=true,
    urlcolor = {darkred},
    linkcolor = {darkred},
    citecolor = {darkred},
    linkcolor = {darkred},
    linktoc=all
}
\usepackage{amssymb}
\usepackage{framed}
\usepackage{amsthm}
\usepackage{graphicx}
\usepackage[all]{xy}
\usepackage{tikz}
\usepackage{tikz-cd}
\usepackage[inline]{enumitem}
\usepackage{mathtools}
\usepackage{etoolbox}
\usepackage{extpfeil}
\usepackage{soul}
\usepackage{authblk}
\usepackage[margin=0.8in]{geometry}
\usepackage{mathtools}

\usetikzlibrary{matrix,arrows,decorations.pathmorphing,decorations.pathreplacing,patterns,shapes.geometric}

\theoremstyle{definition}
\newtheorem{theorem}{Theorem}[section]

\newtheorem{definition}[theorem]{Definition}

\newcommand{\Int}{\mathrm{Int}}

\newcommand{\R}{\mathbb{R}}

\newcommand{\N}{\mathbb{N}}

\newcommand{\Sp}{\mathbb{S}}

\newcommand{\RNum}[1]{\uppercase\expandafter{\romannumeral #1\relax}}

\usepackage[ruled,vlined, linesnumbered]{algorithm2e}

\title{Equichordal Points of Convex Bodies} 

\author[1,*]{Leo Jang}
\author[2,*]{Donghan Kim}

\affil[1]{Department of Computer Sciences, POSTECH, South Korea\thanks{\texttt{leo630@postech.ac.kr} }}
\affil[2]{Department of Mathematical Sciences, KAIST, South Korea\thanks{\texttt{patrick6231@kaist.ac.kr} }}

\begin{document}

\maketitle

\begin{abstract}
The equichordal point problem is a classical question in geometry, asking whether there exist multiple equichordal points within a single convex body. An equichordal point is defined as a point through which all chords of the convex body have the same length. This problem, initially posed by Fujiwara and further investigated by Blaschke, Rothe, and Weitzenböck, has remained an intriguing challenge, particularly in higher dimensions. In this paper, we rigorously prove the nonexistence of multiple equichordal points in $n$-dimensional convex bodies for $n \geq 2$. By utilizing topological tools such as the Borsuk-Ulam theorem and analyzing the properties of continuous functions and mappings on convex bodies, we resolve this long-standing question.
\end{abstract}

\noindent\rule{\textwidth}{0.4pt}
\vspace{-1em} % Adjust spacing as needed
\footnotetext[1]{Equal contribution}

\section{Introduction}
\textbf{Problem.} Does there exist (for $n \geq 2$) an $n$-dimensional convex body that possesses two equichordal points?\cite{klee1960unsolved}

The concept of equichordal points originates from classical geometry and the study of convex bodies. For a set $ C $ that is star-shaped with respect to an interior point $ p $—meaning $ C $ contains every line segment connecting $ p $ to any other point in $ C $—the point $ p $ is called an equichordal point if all chords of $ C $ passing through $ p $ have the same length. A straightforward example is the center of a spherical region, which serves as an equichordal point due to the inherent symmetry of the sphere.

The equichordal point problem for plane convex bodies was first posed by Fujiwara~\cite{fujiwara1916mittelkurve}, who conjectured the existence of multiple equichordal points in certain convex shapes. Independently, Blaschke, Rothe, and Weitzenböck~\cite{blaschke1917} extended this question to more general geometric settings. However, despite substantial progress, the existence of multiple equichordal points in higher-dimensional convex bodies remained an open question.

In this paper, we address the nonexistence of multiple equichordal points in $ n $-dimensional convex bodies ($ n \geq 2 $). By employing tools from topology, such as the Borsuk-Ulam theorem, and leveraging the properties of continuous mappings, we establish that such configurations are impossible.

\section{Preliminaries}\label{sec:preliminaries}
\noindent\textbf{Diameter of a Set.}
The diameter of a set of points in a metric space is the largest distance between points in the set. If $S$ is a set of points with metric $d$, the diameter is
\[
\text{diam}(S)= \sup\limits_{x, y \in S} d(x, y).
\]

\noindent\textbf{Convex Bodies.}
A convex body in $n$-dimensional Euclidean space $\R^n$ is a compact convex set with non-empty interior. 

\noindent\textbf{Chord and Equichordal Point.}
A chord of a circle is a straight line segment whose endpoints both lie on a circular arc. If a chord were to be extended infinitely on both directions into a line, the object is a secant line.
In geometry, an equichordal point is a point defined relative to a convex plane curve such that all chords passing through the point are equal in length. 

\noindent\textbf{Borsuk-Ulam Theorem.}\cite{matouvsek2003using}
The Borsuk–Ulam theorem states that every continuous function from an $n$-sphere into Euclidean $n$-space maps some pair of antipodal points to the same point. Here, two points on a sphere are called antipodal if they are in exactly opposite directions from the sphere's center.
\begin{theorem}[Borsuk-Ulam]\label{lem:Borsuk_Ulam}
    The are no nonconstant antipodal continuous map $f:\Sp^n \to \R^k$ for every $n,k \in \N$ with $k\leq n$.
\end{theorem}

\noindent\textbf{Uniform Metric and Spaces of Continuous Functions.}\cite{munkrestopology}
\begin{definition}
Let $(Y, d)$ be a metric space. Define the \emph{bounded metric} $\bar{d}$ on $Y$ by
\[
\bar{d}(a, b) = \min\{d(a, b), 1\}.
\]
If $\mathbf{x} = (x_\alpha)_{\alpha \in J}$ and $\mathbf{y} = (y_\alpha)_{\alpha \in J}$ are elements of the Cartesian product $Y^J$, the \emph{uniform metric} $\bar{\rho}$ on $Y^J$ is given by
\[
\bar{\rho}(\mathbf{x}, \mathbf{y}) = \sup\{\bar{d}(x_\alpha, y_\alpha) \ | \ \alpha \in J\}.
\]
\end{definition}

For functions $f, g : J \to Y$, the uniform metric $\bar{\rho}$ takes the form
\[
\bar{\rho}(f, g) = \sup\{\bar{d}(f(\alpha), g(\alpha)) \ | \ \alpha \in J\}.
\]

\begin{theorem}
If $(Y, d)$ is a complete metric space, then the product space $Y^J$ is also complete under the uniform metric $\bar{\rho}$.
\end{theorem}

Now consider the subset $\mathcal{C}(X, Y) \subseteq Y^X$ consisting of all \emph{continuous functions} $f : X \to Y$. If $Y$ is complete under the metric $d$, then $\mathcal{C}(X, Y)$ is also complete under the uniform metric $\bar{\rho}$. Similarly, the set $\mathcal{B}(X, Y)$ of \emph{bounded functions} $f : X \to Y$ (where a function $f$ is bounded if $f(X)$ is a bounded subset of $Y$) is complete under the same metric.

\begin{theorem}
Let $X$ be a topological space and $(Y, d)$ a metric space. The set $\mathcal{C}(X, Y)$ of continuous functions and the set $\mathcal{B}(X, Y)$ of bounded functions are both closed subsets of $Y^X$ under the uniform metric. Consequently, if $Y$ is complete under $d$, then $\mathcal{C}(X, Y)$ and $\mathcal{B}(X, Y)$ are complete under $\bar{\rho}$.
\end{theorem}

Given a sequence of functions $(f_n)$ in $\mathcal{C}(X, Y)$ converging to a function $f$ under the uniform metric, for all $x \in X$ and $n \geq N$, we have
\[
\bar{d}(f_n(x), f(x)) \leq \bar{\rho}(f_n, f) < \epsilon,
\]
which implies uniform convergence of $(f_n)$ to $f$.

\begin{definition}
For a metric space $(Y, d)$, another metric on $\mathcal{B}(X, Y)$ is defined as
\[
\rho(f, g) = \sup\{d(f(x), g(x)) \ | \ x \in X\}.
\]
This metric, known as the \emph{sup metric}, is well-defined since the union $f(X) \cup g(X)$ is bounded whenever $f$ and $g$ are bounded functions.
\end{definition}

The relationship between the uniform metric $\bar{\rho}$ and the sup metric $\rho$ is straightforward. For $f, g \in \mathcal{B}(X, Y)$:
\[
\bar{\rho}(f, g) = \min\{\rho(f, g), 1\}.
\]

If $\rho(f, g) > 1$, then there exists $x_0 \in X$ such that $d(f(x_0), g(x_0)) > 1$. Hence, $\bar{d}(f(x_0), g(x_0)) = 1$, and $\bar{\rho}(f, g) = 1$. Conversely, if $\rho(f, g) \leq 1$, then $\bar{d}(f(x), g(x)) = d(f(x), g(x)) \leq 1$ for all $x \in X$, so $\bar{\rho}(f, g) = \rho(f, g)$. Therefore, on $\mathcal{B}(X, Y)$, the uniform metric $\bar{\rho}$ coincides with the bounded version of the sup metric $\rho$.

If $X$ is compact, every continuous function $f : X \to Y$ is bounded, so the sup metric is defined on $\mathcal{C}(X, Y)$. If $Y$ is complete under $d$, then $\mathcal{C}(X, Y)$ is complete under both the uniform metric $\bar{\rho}$ and the sup metric $\rho$. In practice, the sup metric is often preferred in this context.

\section{Nonexistence of Multiple Equichordal Points in Convex Bodies}
In this section, we prove the nonexistence of multiple equichordal points in a convex body $ X $ for $ n \geq 2 $.

Let $ X $ be a convex body in $n$-dimensional Euclidean space $\R^n$ for $n \geq 2$. Pick any point $ x \in \Int X $. Since $x$ is an interior point of $X$, there exists $\epsilon > 0$ such that $B_d(x, \epsilon) \subset X$. 
Thus, we can consider a unit vector $a$ on $ S^{n-1} $, the unit sphere in $n$-dimensional space. 
Let $d$ be a fixed length of the half-chord starting at 
$x$ in the given direction $a$. Thus, we have
\[
d = \sup\limits_{t>0} \{d(x, x+at) \ | \ x+at \in X \}.
\]
Since $d(x, x+at)$ is a continuous function of $\R_{>0}$ and $B$ is a closed set, there exists a unique $t_0>0$ such that $d=d(x, x+at_0)$ and $a_0:=x+at_0$ lies on $B$ such that the chord in the given direction starting at $x$ has a fixed length $d$, where $d$ is the distance from the given interior point to the point $a_0$.

Define the function
\[
\varphi: \Int X \to \mathcal{C}(\Sp^{n-1}, \R)
\]
with
\[
\varphi(x)(a) = d(x, a_0),
\]
where $ \varphi(x) $ is a function from $ \Sp^{n-1} $ to $ \R $. We will show that $ \varphi $ is a continuous injective map.

\subsection{Continuity of $ \varphi(x) $ for each $ x \in \Int X$}
For a fixed $ x \in X $ and any $ a, b \in \Sp^{n-1} $, we have
\[
|\varphi(x)(a) - \varphi(x)(b)| = |d(x, a_0) - d(x, b_0)| \leq d(a_0, b_0).
\]
Since,
\[
d(a_0, b_0) \leq  \text{diam}(X)d(a, b),
\]
where the inequality holds because the maximum distance between two points on $ X $ is less than the geodesic distance. Since $d(x, a) \leq \text{diam}(X)$ and $d(x, b) \leq \text{diam}(X)$, the geodesic distance between $a, b$ is less than $\text{diam}(X)d(a, b)$. Thus, we have $|\varphi(x)(a) - \varphi(x)(b)| \leq \text{diam}(X)d(a, b)$. Since $X$ is a compact set, $\text{diam}(X)$ is finite. Hence, $ \varphi(x) $ is a $1$-Lipschitz map and therefore a continuous function.

\subsection{Injectivity of $ \varphi $}
Suppose $ x, y \in \Int X $ are distinct points. Let $ a \in \Sp^{n-1} $ be the unit vector in the direction from $ x $ to $ y $. i.e., $a=\frac{y-x}{||y-x||}$. Then, for the chord passing through $ x $ and $ y $, we have
\[
d(x, a_0) = d(y, a_0) + d(x, y).
\]
This implies that $ \varphi(x)(a_0) \neq \varphi(y)(a_0) $, so $ \varphi(x) \neq \varphi(y) $. Therefore, $ \varphi $ is injective.

\subsection{Continuity of $ \varphi $}
To show $ \varphi $ is continuous, let $ x, y \in \Int X $ and consider any $ a \in \Sp^{n-1} $. Then,
\[
|\varphi(x)(a) - \varphi(y)(a)| = |d(x, a_0) - d(y, a_0)| \leq d(x, y).
\]
Since the inequality holds for all $ a \in \Sp^{n-1} $, we have
\[
\rho(\varphi(x), \varphi(y)) \leq d(x, y),
\]
where $ \rho $ is the uniform metric on $ \mathcal{C}(\Sp^{n-1}, \R) $. Thus, $ \varphi $ is a $ 1 $-Lipschitz map and therefore a continuous function.

\subsection{Nonexistence of Multiple Equichordal Points}
Suppose there exist two distinct equichordal points $ x, y \in \Int X $. Let $ a \in \Sp^{n-1} $ be the unit vector from $ x $ to $ y $, and let $ b = -a $. The length of the chord passing through $ x $ is given by
\[
d(x, a_0) + d(x, b_0)
\]
and for $ y $, it is
\[
d(y, a_0) + d(y, b_0) = d(x, a_0) - d(x, y) + d(x, b_0) + d(x, y) = d(x, a_0) + d(x, b_0).
\]
Thus, the chord lengths at both $ x $ and $ y $ are equal, with a common length $ r $. For every $ s \in \Sp^{n-1} $, we have
\[
\varphi(x)(s) + \varphi(x)(-s) = r.
\]
Define the function
\[
g: \Sp^{n-1} \to \R, \quad g(s) = \varphi(x)(s) - \varphi(y)(s).
\]
Since $ \varphi(x) $ and $ \varphi(y) $ are continuous, $ g $ is also continuous. Moreover,
\begin{align*}
    g(-s) &= \varphi(x)(-s) - \varphi(y)(-s) \\
    &= (r-\varphi(x)(s))-(r-\varphi(y)(s)) \\
    &= -(\varphi(x)(s)-\varphi(y)(s)) \\
    &= -g(s).
\end{align*}
Thus, $ g $ is an antipodal continuous map on $ \Sp^{n-1} $. By the Borsuk-Ulam theorem (Theorem~\ref{lem:Borsuk_Ulam}), $ g(s) = 0 $ for all $ s \in \Sp^{n-1} $, implying $ \varphi(x) = \varphi(y) $.

Since $ \varphi $ is injective, this implies $ x = y $, which is a contradiction. Therefore, there cannot exist multiple equichordal points in $ X $.

\qed

\bibliographystyle{plainurl}
% \bibliography{bib.bib} 

\begin{thebibliography}{1}

\bibitem{blaschke1917}
W.~Blaschke, H.~Rothe, and R.~Weitzenbock.
\newblock Aufgabe 552.
\newblock {\em Archiv der Math. u. Physik}, 27:82, 1917.

\bibitem{fujiwara1916mittelkurve}
Matsusaburo Fujiwara.
\newblock {\"U}ber die mittelkurve zweier geschiossenen konvexen kurven in bezug auf einen punkt.
\newblock {\em Tohoku Mathematical Journal, First Series}, 10:99--103, 1916.

\bibitem{klee1960unsolved}
Victor Klee.
\newblock Some unsolved problems in geometry.
\newblock {\em American Mathematical Monthly}, 1960.

\bibitem{matouvsek2003using}
Ji{\v{r}}{\'\i} Matou{\v{s}}ek, Anders Bj{\"o}rner, G{\"u}nter~M Ziegler, et~al.
\newblock {\em Using the Borsuk-Ulam theorem: lectures on topological methods in combinatorics and geometry}, volume 2003.
\newblock Springer, 2003.

\bibitem{munkrestopology}
James~R. Munkres.
\newblock {\em Topology}.
\newblock Prentice Hall, 2000.

\end{thebibliography}

\end{document}